\documentclass[11pt,oneside]{amsart}

\usepackage{latexsym,amsmath}
\usepackage{amsmath,amsthm}
\usepackage{amsfonts}
\usepackage[psamsfonts]{amssymb}
\usepackage{graphicx}
\usepackage{pstricks}

\theoremstyle{plain}

\newtheorem{proposition}{Proposition}

\theoremstyle{definition} 

\theoremstyle{definition} 

\theoremstyle{remark} 

\theoremstyle{remark} 

\newtheorem*{remark*}{Remark}
\numberwithin{equation}{section}

\renewcommand{\u}{\,%
\psset{unit=.7pt}
\begin{pspicture}(0,0)(9,10)
\psline[linewidth=.3pt,origin={0,1}]{-}(0,0)(2,5)(7,5)(9,10)
\end{pspicture}
}

\renewcommand{\d}{\,%
\psset{unit=.7pt}
\begin{pspicture}(0,0)(9,10)
\psline[linewidth=.3pt,origin={0,1}]{-}(0,10)(2,5)(7,5)(9,0)
\end{pspicture}
}

\newcommand{\ud}{\,
\psset{unit=.7pt}
\begin{pspicture}(0,0)(16,10)
\psline[linewidth=.3pt,origin={0,1}]{-}(0,0)(4,10)(12,10)(16,0)
\end{pspicture}
}

\newcommand{\du}{\,
\psset{unit=.7pt}
\begin{pspicture}(0,0)(16,10)
\psline[linewidth=.3pt,origin={0,1}]{-}(0,10)(4,0)(12,0)(16,10)
\end{pspicture}
}

\newcommand{\uudd}{\,
\psset{unit=.7pt}
\begin{pspicture}(0,0)(23,10)
\psline[linewidth=.3pt,origin={0,1}]{-}(0,0)(2,5)(7,5)(9,10)(14,10)(16,5)(21,5)(23,0)
\end{pspicture}
}

\newcommand{\dduu}{\,
\psset{unit=.7pt}
\begin{pspicture}(0,0)(23,10)
\psline[linewidth=.3pt,origin={0,1}]{-}(0,10)(2,5)(7,5)(9,0)(14,0)(16,5)(21,5)(23,10)
\end{pspicture}
}

\renewcommand{\le}{\leqslant}

\newcommand{\tr}{\tilde\rho}

\newcommand{\R}{\mathbb{R}}
\newcommand{\Z}{\mathbb{Z}}

\newcommand{\sign}{\operatorname{sign}}

\newcommand{\ru}{\rule[-6pt]{0pt}{17pt}}

\title[L'Hospital rules for monotonicity]{``Non-strict'' l'Hospital-Type Rules for Monotonicity: Intervals of Constancy}
\author{Iosif Pinelis}
\address{Department of Mathematical Sciences, Michigan Technological University, Hough\-ton, Michigan 49931}
\email{ipinelis@mtu.edu}

\begin{document}

\today\quad\jobname.tex


\maketitle


Let $f$ and $g$ be 
differentiable functions
defined on the interval $(a,b)$, where \break 
$-\infty\le a<b\le\infty$, and let 
$$r:=\frac fg\quad\text{and}\quad\rho:=\frac{f'}{g'}.$$  
It is assumed throughout 
that $g$ and $g'$ do not take on the zero value 
anywhere on $(a,b)$.
The function $\rho$ may be referred to as \emph{a derivative ratio} for the ``original'' ratio $r$. 
In \cite{3}, general ``rules" for monotonicity patterns, resembling the usual l'Hospital rules for limits, were given.
In particular, according to \cite[Proposition~1.9 and Remark~1.14]{3}, 
one has the dependence of the monotonicity pattern of $r$ \big(on $(a,b)$\big) on that of 
$\rho$
(and also on the sign of $gg'$) as given by Table~\ref{t:1}. 
The vertical double line in the table separates the conditions (on the left) from the corresponding conclusions (on the right). 
\begin{table}[h] 
\begin{center}
\begin{tabular}{c|c||c}
$\rho$ \ru	& $gg'$	& $r$ \\ \hline\hline
$\u$ \ru		& $>0$	& $\dduu$ \\ \hline
$\d$ \ru		& $>0$	& $\uudd$ \\ \hline
$\u$ \ru		& $<0$	& $\uudd$ \\ \hline
$\d$ \ru		& $<0$	& $\dduu$ \\ \hline 
\end{tabular}
\vspace*{10pt}
\caption[smallcaption]{``Non-strict'' general rules for monotonicity.}
\label{t:1}
\end{center}
\end{table}

Here, for instance, $r\dduu$ means that there is some $c\in[a,b]$ such that $r\d$ (that is, $r$ is non-increasing) on $(a,c)$ and 
$r\u$ ($r$ is non-decreasing) on $(c,b)$; 
in particular, if $c=a$ then $r\dduu$ simply means that $r\u$ on the entire interval $(a,b)$; and if $c=b$ then $r\dduu$ means that $r\d$ on $(a,b)$.
Thus, if one also knows whether $r\u$ or $r\d$ in a right neighborhood of $a$ and in a left neighborhood of $b$, then Table~\ref{t:1} uniquely determines the ``non-strict'' monotonicity pattern of $r$. 
\big(The ``strict'' counterparts of these rules, with terms ``increasing'' and ``decreasing'' in place of ``non-decreasing'' and ``non-increasing'' respectively, also hold, according to the same Proposition~1.9 of \cite{3}.\big)

Clearly, the stated l'Hospital-type rules for monotonicity patterns are helpful wherever the l'Hospital rules for limits are so, and even beyond that, because these monotonicity rules do not require that both $f$ and $g$ (or either of them) tend to 0 or $\infty$ at any point. (Special rules for monotonicity, which do require that both $f$ and $g$ vanish at an endpoint of $(a,b)$, were given, in different forms and with different proofs, in \cite{grom1,pin91,avv93,2,borwein}.)

Thus, it should not be surprising that a wide variety of applications of the l'Hospital-type rules for monotonicity patterns were given: in areas of
analytic inequalities \cite{avv05,2,3,monthly}, 
approximation theory \cite{4}, 
differential geometry \cite{chavel,grom1,grom2,saccheri},
information theory \cite{2,3}, 
(quasi)conformal mappings \cite{aqvv,avv93,avv97,avv01}, 
probability and statistics \cite{pin91,3,4,5,binom,normal,edelman,houdre}, 
etc. 

These rules for monotonicity could be helpful when $f'$ or $g'$ can be expressed simpler than or similarly to $f$ or $g$, respectively. 
Such functions $f$ and $g$ are essentially the same as the functions that could be taken to play the role of $u$ in the integration-by-parts formula $\int u\,dv=uv-\int v\,du$;
this class of functions includes algebraic, exponential, trigonometric, logarithmic, inverse trigonometric and inverse hyperbolic functions, and as well as non-elementary ``anti-derivative" functions of the form $x\mapsto c+\int_a^x h(u)\,du$ or $x\mapsto c+\int_x^b h(u)\,du$. 

``Discrete'' analogues, for $f$ and $g$ defined on $\Z$, of the l'Hospital-type rules for monotonicity are available as well \cite{discrete}. 

In this paper, we describe different facets of the relation of the (maximal) interval(s) of constancy of the original ratio $r$ with those of the derivative ratio~$\rho$.

$\P$ In what follows, let us always assume that 
$\rho$ is (not necessarily strictly) monotonic (that is, $\u$ or $\d$) on $(a,b)$.

Let us say that an interval $I\subseteq(a,b)$ is an \emph{interval of constancy (i.c.)} of a function $h\colon(a,b)\to\R$ if $I$ is of nonzero length and $h$ is constant on $I$. If an i.c.\ $I$ is not contained in any other i.c., let us say that $I$ is a \emph{maximal i.c.\ (m.i.c.)} It is easy to see that any i.c.\ is contained in a unique m.i.c.\ (which is simply the union of all i.c.'s containing the given i.c.).

It is easy to see that every i.c.\ of $r$ is an i.c.\ of $\rho$. 
One might think that, if $\rho$ has more than one m.i.c., then this can also be the case for the original ratio $r$. It may therefore be surprising that the opposite is true, and even in the following strong sense. 

\begin{proposition}\label{prop:strong}
The rules given by Table~\ref{t:1} can be strengthened as shown in Table~\ref{t:2}.

  \begin{table}[h] 
\begin{center}
\begin{tabular}{c|c||c}
$\rho$ \ru	& $gg'$	& $r$ \\ \hline\hline
$\u$ \ru		& $>0$	& $\du$ \\ \hline
$\d$ \ru		& $>0$	& $\ud$ \\ \hline
$\u$ \ru		& $<0$	& $\ud$ \\ \hline
$\d$ \ru		& $<0$	& $\du$ \\ \hline 
\end{tabular}
\vspace*{10pt}
\caption[smallcaption]{Improved ``non-strict'' general rules for monotonicity.}
\label{t:2}
\end{center}
\end{table}

\vspace*{-.5cm}
\noindent
Here, for instance, $r\du$ means that there is a subinterval $[c,d]\subseteq[a,b]$ (possibly of length $0$) such that $r'<0$ on $(a,c)$, $r$ is constant on $(c,d)$, and 
$r'>0$ on $(d,b)$.
\end{proposition}

Why is this proposition true? The key notion here is that of the function 
$$
\tr:=r'\,\frac{g^2}{|g'|}
,	
$$
introduced in \cite{3} and further studied in \cite{borwein}. 
The key lemma concerning $\tr$ \cite[Lemma~1 and Remark~4]{borwein} states, as presented in Table~\ref{t:3} here, that the monotonicity pattern of $\tr$ is the same as that of $\rho$ if $gg'>0$, and opposite to the pattern of $\rho$ if $gg'<0$. 
\begin{table}[h] 
\begin{center}
\begin{tabular}{c|c||c}
$\rho$ \ru	& $gg'$	& $\tr$ \\ \hline\hline
$\u$ \ru		& $>0$	& $\u$ \\ \hline
$\d$ \ru		& $>0$	& $\d$ \\ \hline
$\u$ \ru		& $<0$	& $\d$ \\ \hline
$\d$ \ru		& $<0$	& $\u$ \\ \hline 
\end{tabular}
\vspace*{10pt}

\caption{The monotonicity patterns of $\rho$ and $\tr$ mirror each other}
\label{t:3}
\end{center}
\end{table}


From this relation between $\rho$ and $\tr$, the rules given by Table~\ref{t:1} can be easily deduced, since 
$$\sign(r')=\sign\tr.$$

A simple but important observation is that the derivative ratio $\rho$ and its counterpart $\tr$ are continuous functions \cite[Remark~4]{borwein}. Since the ratio $r$ is differentiable, it is continuous as well. Therefore, any m.i.c.\ $I$ of $r$ or $\rho$ or $\tr$ is closed (as a set) in $(a,b)$; that is, $I$ has the form $[c,d]$ or $(a,c]$ or $[d,b)$ or $(a,b)$, for some $c$ and $d$ such that $a<c<d<b$. 
Moreover, it is seen from Table~\ref{t:3} that the m.i.c.'s of $\tr$ are the same as those of $\rho$, because any real function $h$ is constant on an interval $I$ if and only if $h$ is both non-decreasing and non-increasing on $I$.

Proposition~\ref{prop:strong} now follows easily.
Indeed, it suffices to consider only the first line of Table~\ref{t:2} \big(since the other three lines can then be obtained by ``vertical'' reflection $f\leftrightarrow-f$ and/or ``horizontal'' reflection $x\leftrightarrow-x$\big).
So, assume that $\rho\u$ on $(a,b)$ and $gg'>0$. 
Then $\tr\u$ on $(a,b)$. If $\tr>0$ and hence $r'>0$ on the entire interval $(a,b)$, let $c:=d:=a$, to obtain the conclusion that $r\du$ on $(a,b)$. If $\tr<0$ and hence $r'<0$ on $(a,b)$, let $c:=d:=b$. 
It remains to consider the case when the sign of $\tr$ takes on at least two different values (of the set $\{-1,0,1\}$ of all sign values). 
Then, since the function $\tr$ is non-decreasing and continuous on $(a,b)$,
the level-$0$ set $\ell_0(\tr):=\{u\in(a,b)\colon\tr(u)=0\}$ of $\tr$ must be a non-empty interval \big(which in fact must be an m.i.c.\ of $\tr$ and hence a set closed in $(a,b)$\big); in this case, take $c$ and $d$ to be the left and right endpoints, respectively, of the interval $\ell_0(\tr)$ (at that, it is possible that $c=a$ and/or $d=b$). Then $\tr<0$ and hence $r'<0$ on $(a,c)$; $\tr=0$ and hence $r'=0$ and $r=\text{const}$ on $(c,d)$; and $\tr>0$ and hence $r'>0$ on $(d,b)$. 

By Proposition~\ref{prop:strong}, $r$ can have no more than one m.i.c. On the other hand, one has 

\begin{proposition}\label{prop:2}
If $r$ has an m.i.c.\ $I$, then $I$ must be an m.i.c.\ of $\rho$ and $\tr$ as well. 
\end{proposition}

Indeed, suppose that $I$ is the (necessarily unique) m.i.c.\ of $r$, so that $\frac fg=r=K$ on $I$ for some constant $K$. Then obviously $\rho=K$ and $\tr=0$ on $I$, so that $I$ is an i.c.\ of $\rho$ and $\tr$. Let then $J$ be the unique m.i.c.\ of $\rho$ such that $J\supseteq I$, whence $\frac{f'}{g'}=\rho=K_1$ on $J$ for some constant $K_1$, and so, $f=K_1 g+C$ and $r=K_1+\frac Cg$ on $J$, and hence on $I$, for some constant $C$. But $r$ is constant on the nonzero-length interval $I$, while $g$ is not constant on $I$ (because $g'(x)\ne0$ for any $x\in(a,b)$). It follows that $C=0$ and thus $r=K_1$ on $J$. Finally, since $I$ is an m.i.c.\ of $r$ and $J\supseteq I$, one concludes that $J=I$, and so, $I$ is an m.i.c.\ of $\rho$ and hence of $\tr$.

We complete the description of the relation between the m.i.c.'s of $r$ and $\rho$ by observing that any one m.i.c.\ $I$ of a given derivative ratio $\rho$ is the m.i.c.\ of an appropriately constructed original ratio $r$ (which must, in view of Proposition~\ref{prop:strong}, depend on the choice of $I$):

\begin{proposition}\label{prop:4}
For 
\begin{trivlist}
\item
\ $\bullet$\ any differentiable function $g\colon(a,b)\to\R$ such that  $gg'(x)\ne0$ for each $x\in(a,b)$, 
\item
\ $\bullet$\ any (not necessarily strictly) monotonic continuous function $\rho\colon(a,b)\to\R$, and 
\item
\ $\bullet$\ any m.i.c.\ $I$ of $\rho$
\end{trivlist}
there exists a differentiable function $f\colon(a,b)\to\R$ such that $\dfrac{f'}{g'}=\rho$ and the only m.i.c.\ of $r:=\dfrac fg$ is $I$. 
\end{proposition}

Indeed, let $g$, $\rho$, and $I$ satisfy the conditions listed in Proposition~\ref{prop:4}, so that $\rho=K$ on $I$ for some constant $K$.  
Note that the condition on $g$ implies that either $g'>0$ on the entire interval $(a,b)$ or $g'<0$ on $(a,b)$ (see e.g.\ \cite[Remark~3]{borwein}), so that $g$ is monotonic and hence of locally bounded variation on $(a,b)$.
Take any point $z$ in the interval $I$ (which is an i.c.\ and hence non-empty) and define $f$ by the formula
$$f(x):=Kg(z)+\int_z^x\,\rho(u)\,d\,g(u)$$
for all $x$ in $(a,b)$, where the integral may be understood in the 
Riemann-Stieltjes sense, with the convention that $\int_z^x:=-\int_x^z$ if $x<z$. 
Because $\rho$ is continuous and $g$ is differentiable, it follows that for the so defined function $f$ one has $\frac{f'}{g'}=\rho$; 
moreover, $f=Kg$ on $I$, so that $I$ is an i.c.\ of $r$. 
But any i.c.\ of $r$ is also an i.c.\ of $\rho$, and $I$ was assumed to be an m.i.c.\ of $\rho$. It follows that $I$ is \emph{an} m.i.c. (and hence \emph{the only} m.i.c.) of $r$.

Let us summarize our findings: (i) the set of all m.i.c.'s of $\tr$ is the same as that of $\rho$; (ii) $r$ can have at most one m.i.c., and its m.i.c.\ 
must also be an m.i.c.\ of $\rho$ and thus of $\tr$; moreover, then the m.i.c. of $r$ is the level-$0$ set of $\tr$; (iii) any one m.i.c. of a given derivative ratio $\rho$ is the m.i.c.\ of an appropriately constructed original ratio~$r$.

\end{document}